\documentclass[preprint,12pt]{elsarticle}

\usepackage{amssymb}
\usepackage{amsthm}
\usepackage{amsmath}

\usepackage{algorithm}
\usepackage{algorithmic}

\newtheorem{thm}{Theorem}
\newtheorem{lem}[thm]{Lemma}
\newtheorem{prop}[thm]{Proposition}
\newtheorem{cor}[thm]{Corollary}
\newproof{pf}{Proof}

\begin{document}

\begin{frontmatter}



\title{Convex $p$-partitions of bipartite graphs}


\author[ICI]{Luciano~N.~Grippo}\ead{lgrippo@ungs.edu.ar}
\author[DIM,CMM]{Mart{\'i}n~Matamala}\ead{mmatamal@dim.uchile.cl}
\author[ICI]{Mart{\'i}n~D.~Safe}\ead{msafe@ungs.edu.ar}
\author[CMM]{Maya~J.~Stein}\ead{mstein@dim.uchile.cl}

\address[ICI]{Instituto de Ciencias, Universidad Nacional de General Sarmiento, Los Polvorines, Buenos Aires, Argentina}
\address[DIM]{Departamento de Ingenier{\'i}a Matem{\'a}tica, Universidad de Chile, Santiago, Chile}
\address[CMM]{Centro de Modelamiento Matem{\'a}tico (CNRS-UMI 2807), Universidad de Chile, Santiago, Chile}

\begin{abstract}
A set of vertices $X$ of a graph $G$ is \emph{convex} if no shortest path
between two vertices in $X$ contains a vertex outside $X$.
We prove that for fixed $p\geq 1$,  all
partitions  of the vertex set of a bipartite graph into $p$ convex sets can be found in polynomial time.
\end{abstract}


\end{frontmatter}


\section{Introduction}

Given a graph $G=(V,E)$, a set $X$ of vertices is called \emph{convex}
if $G[X]$, the graph induced by $X$, contains all shortest paths between any two of its vertices. All graphs here are undirected and simple. The notion probably first appeared in~\cite{FeldH}, see also~\cite{HN1981}, and later became also known as  geodesic convexity, or d-convexity, in order to distinguish it from  different notions of convexity in graphs and other combinatorial structures (see~\cite{D1987} for an early overview). The book~\cite{Pelayo} gives an up-to-date survey of results on convexity in graphs.

One of the approaches to convexity in graphs comes
from the viewpoint of computational complexity. Clearly, by computing the distances between all pairs, one can decide in polynomial time if a given set of vertices is convex. To determine the size of a largest convex set not covering the whole graph, however, is an NP-complete problem, even for bipartite graphs, albeit linear for cographs~\cite{DPRS2012}. The same phenomenon occurs (NP-completeness even for bipartite graphs, but linearity for cographs) if we wish to determine related invariants such as the hull number and the geodetic number of a graph~\cite{A,D, DPS2008}.

We focus here on the notion of a \emph{convex $p$-partition} of a graph, that is, a partition of the vertex set  into $p$ convex sets.
For instance, any graph on~$n$ vertices containing a matching of size $m$ has a convex $(n-m)$-partition, and trivially, any graph has a convex $1$-partition.
Deciding whether a graph has a convex $p$-partition, for fixed $p \geq 2$, is NP-complete for arbitrary graphs, and linear time solvable for cographs~\cite{addsDM2011}. Also, any connected chordal graph has at least one convex $p$-partition for each $p\geq 1$~\cite{addsDM2011}.

In view of the above described panorama, it was conjectured in~\cite{Pelayo} that also for bipartite graphs, it should be NP-complete to decide whether they have a convex $p$-partition.
We show that, for any fixed $p\geq 1$,  this is not the case. More precisely, we prove that for  $p\geq 1$, all convex $p$-partitions of a bipartite graph can be enumerated in polynomial time.
This extends a recent result of Glantz and Mayerhenke \cite{GM2013}, who prove the same for the case $p=2$. They also showed that all convex $2$-partitions of a planar graph can be found in polynomial time.

\section{Bipartite graphs with convex $p$-partitions}

We start by reproving the result for bipartite graphs from \cite{GM2013} in a slightly different way. At the same time, this will serve as a base for the general case. We denote the distance between two vertices $u$ and $v$ in a graph $G$ by $d_G(u,v)$, defined as the length of a shortest path between $u$ and $v$. It is known that for a given $u$, the set of all distances $d(u,v)$, for $v\in V$, can be computed in linear time (\cite{Thorup1999}).

 \begin{lem}\label{l:frontierdistances}
  Given a convex set $C$ in a connected bipartite graph $G$,
  and an edge $uv$ with $u\in C$, $v\notin C$ we have that $d_G(u',u)<d_G(u',v)$,
  for each $u'\in C$.
 \end{lem}

 \begin{pf}
 Suppose otherwise. Observe that since $G$ is bipartite, $d_G(u',u)\neq d_G(u',v)$, and thus we may assume $d_G(u',u)>d_G(u',v)$.  Then there is a shortest path $P$ from $u'$ to $v$ not containing $u$. Extending $P$ to $u$ through the edge $vu$, gives a shortest path from $u'$ to $u$, a contradiction, as $u$  and $u'$ lie in the convex  set $C$, but
 $v\notin C$.\qed
 \end{pf}

 Let $e=uv$ be an edge of $G$ and denote by $X_{uv}$ the set
of vertices that are closer to $u$ than to $v$. If $G$ is a connected bipartite graph, then	
$V$ is the disjoint union $X_{uv}\cup X_{vu}$. From Lemma
\ref{l:frontierdistances} we get the following corollaries.

\begin{cor}\label{c:edgeinconvexcut}
Let $uv$ be an edge of a connected bipartite graph $G$.
If $C$ is a convex set containing $u$ and not containing
$v$, then $C\subseteq X_{uv}$.
\end{cor}

\begin{cor}\label{c:edgebetweentwoconvexs}
Let $G=(V,E)$ be a connected bipartite graph, with a partition of $V$ into convex sets $X_1, X_2$. Let $uv\in E$, with $u\in X_1$ and $v\in X_2$. Then $X_1\subseteq X_{uv}$ and $X_2\subseteq X_{vu}$ which, as $V=X_{uv}\cup X_{vu}$,  implies that $X_1= X_{uv}$ and $X_2= X_{vu}$.
\end{cor}

From the previous corollary it is direct that there are at most $|E|$ convex
2-partitions and, as a consequence, we can enumerate all
convex $2$-partitions in polynomial time.

\begin{prop}
We can enumerate in polynomial time  all convex 2-partitions of a connected bipartite graph.
\end{prop}

 We now prove that for fixed $p\geq 3$, we can enumerate in polynomial time
all convex $p$-partitions of a connected bipartite graph.
In order to do so,  we extend the idea present in Corollary~\ref{c:edgebetweentwoconvexs}.

 We write $[p]$ for the set $\{1,\ldots,p\}$.
For a set $F$ of edges, let  $V(F)$ denote the set of all endvertices of edges of $F$.

Given a convex $p$-partition $\mathcal X=\{X_1,X_2,\ldots, X_p\}$ of a graph $G=(V,E)$, we call a pair $(F,\phi)$ an \emph{$\mathcal X$-skeleton}, if $F\subseteq E$ and $\phi : V(F)\to [p]$ satisfy the following:
\begin{itemize}
\item all edges of $F$ go between distinct parts of $\mathcal X$;
\item if there is at least one edge in $E$ between $X_i$ and $X_j$, then there is exactly one edge of $F$ between $X_i$ and $X_j$;
\item   $\phi (v)=i$ if and only if  $v\in X_i$.
\end{itemize}
Note that the first two conditions might be equivalently expressed by saying that after contracting the sets $X_i$ and deleting all remaining edges that are not in $F$, we are left with a (simple) graph $H_{(F,\phi)}$ whose edges represent the edges of $G$ that cross the partition. The last condition says $\phi$ assigns the same colour to all vertices of $V(F)$ that become identified in $H_{(F,\phi)}$.

Note that for a connected graph $G$ the second condition implies that $V(F)\cap X_j\neq \emptyset$, for each $j\in [p]$. Then, the third condition implies that $\phi$ is a surjective
function.

Given a set of edges $F$ we say that a function $\phi:V(F)\to [p]$ is
a \emph{$p$-coloring} of $F$  if it is surjective and for each $vw\in F$, $\phi(v)\neq \phi(w)$.

We shall prove that given a graph $G=(V,E)$, $F\subseteq E$
and $\phi$ a $p$-coloring of $F$,
we can decide in linear time  whether $(F,\phi)$ is  the $\mathcal X$-skeleton of a convex $p$-partition $\mathcal X=\{X_1,X_2,\ldots, X_p\}$ of $G$.

To this end, we use the following two criteria which follow from Corollary~\ref{c:edgeinconvexcut} and the definition of a convex set, respectively.

\begin{enumerate}
\item For each $i\in [p]$ and for each edge $vw\in F$ with $\phi(w)=i$,
if a vertex $u\in X_{vw}$ then $u\notin X_{i}$.

\item For each $i\in [p]$, for any three distinct vertices $u,v,w$ with $w\in V(F)$ and $d_G(u,w)=d_G(u,v)+d_G(v,w)$, if $v\notin X_i$ and $w\in X_i$, then $u\notin X_i$.
\end{enumerate}

\begin{algorithm}
\caption{$p$-Is-Skeleton}
\label{al:isskeleton}
\begin{algorithmic}\label{alg}
\REQUIRE A graph $G=(V,E)$ and $(F,\phi)$, $F\subseteq E$, $\phi$ a $p$-coloring of $F$.

\hskip -0.4cm{\bf Return:} $\mathcal X$ -- a convex $p$-partition of $G$ having $(F,\phi)$ as skeleton, if it exists.

\hskip -0.4cm{\bf For all} {$i\in [p]$}

\hskip -0.4cm\ \ \ \  $X^1_i\leftarrow V$;

\hskip -0.4cm\ \ \ \ {\bf For all} {$vw\in F$ with $\phi(w)=i$}

\hskip -0.4cm\ \ \ \ \ \ \ \ \ {\bf For all} {$u\in X_{vw}$}

\hskip -0.4cm\ \ \ \ \ \ \ \ \ \ \ \ \ $X^1_i\leftarrow X^1_i\setminus\{u\}$;

\hskip -0.4cm\ \ \ \ \ \ \ \ \ \ \ \ \ {\bf If} $X^1_i=\emptyset$ {\bf then} {\bf return} False.

\hskip -0.4cm{\bf For all} {$i\in [p]$}

\hskip -0.4cm\ \ \ \  $X^2_i\leftarrow X^1_i$;

\hskip -0.4cm\ \ \ \ {\bf For all} {$w\in V(F)$ with $\phi(w)=i$}

\hskip -0.4cm\ \ \ \ \ \ \ \ \ {\bf For all} {$u\in X^1_i$ s.t. $\exists v\in V\setminus X^1_i$ with $d_G(u,v)+d_G(v,w)=d_G(u,w)$}

\hskip -0.4cm\ \ \ \ \ \ \ \ \ \ \ \ \ $X^2_i\leftarrow X^1_i\setminus\{u\}$;

\hskip -0.4cm\ \ \ \ \ \ \ \ \ \ \ \ \ {\bf If} $X^2_i=\emptyset$ {\bf then} {\bf return} False.

\hskip -0.4cm{\bf For all} {$i\in [p]$}

\hskip -0.4cm\ \ \ \ \ \ \ \ \ \ \ \ \ {\bf If} $X^2_i$ is not convex {\bf then} {\bf return} False.

\hskip -0.4cm{\bf Return:} $\mathcal X=\{X^2_1,\ldots,X^2_p\}$
\end{algorithmic}
\end{algorithm}

The algorithm described in Algorithm \ref{alg} has three steps.
It starts with considering for each part of the convex partition the whole set of vertices.
In a second step, it eliminates from each part $X_i$ those vertices indicated
by the first criterion. For each $vw\in F$ we can compute in linear time
the set $X_{vw}$, and thus,  we can check in constant time whether $u\in X_{vw}$.
Therefore, this part takes linear time.
Finally, in the third step, the algorithm eliminates all vertices indicated by the second criterion from the parts obtained in the previous
step. As before,
for each $w\in V(F)$ we can compute, in linear time, the distance from $w$ to all the vertices, and during the same process,  we can already eliminate the vertices that are as in the second criterion.
Hence, Algorithm \ref{alg} runs in linear time.

The correctness of this algorithm is proved in the next result.
However, in the proof, instead of working with parts $X_i$, we associate to each
vertex the set of indices of the parts to which it belongs. Initially this set is $[p]$,
in the second step we erase for these sets all the indices indicated by the first criterion,
and in the third step we erase from the remaining indices those indicated by the second criterion.

 For each pair of vertices $u$ and $w$ in a graph $G$
 we define $I[u,w]$ as the set of vertices in shortest paths between $u$ and $w$.
 Then, $v\in I[u,w]$ if and only if $d_G(u,w)=d_G(u,v)+d_G(v,w)$.

\begin{thm}\label{t:main}
Let $G=(V,E)$ be a connected bipartite graph, let $F\subseteq E$ and let $\phi : V(F)\to [p]$. If $G$ has a convex $p$-partition with skeleton $(F,\phi)$, then this partition is unique. We can find such partition, or show it does not exist, in polynomial time.
\end{thm}

\begin{pf}
Define lists $L(u)$ for each vertex $u\in V$ by setting $$L(u):=[p]-\{\phi(w)\colon u\in X_{vw}\text{ for some }vw\in F\}.$$ The idea behind the lists $L(u)$ is that they do not contain colours correspon\-ding to  partition sets $u$ cannot belong to, as explained in more detail above. Restricting these lists even more, we define, for each vertex
$u$,
$$L'(u):=L(u)-\{\phi(w)\colon w\in V(F)\mbox{ and }\phi(w)\notin L(v)\text{ for some }v\in I[u,w]\}.$$

We will prove that if $G$ has a convex $p$-partition $\mathcal X=\{X_1,\ldots,X_p\}$ with skeleton $(F,\phi)$,
then, for each $i\in[p]$,
\begin{equation}\label{eq:2}
\text{$L'(u)=\{i\}$ for every $u\in X_i$.}
\end{equation}

We first observe that
\begin{equation}\label{eq:3}
  \mbox{$i\in L(u)$ for every $u\in X_i$}.
\end{equation}
Otherwise, there are $u\in X_i$ and $vw\in F$ such that $\phi(w)=i$ and $u\in X_{vw}$. Hence, $u,w\in X_i$ and $v\in I[u,w]$. Since $X_i$ is convex, $v\in X_i$, contradicting the fact that the edge $vw$ of $F$ must join distinct parts of $\mathcal X$. This contradiction proves~\eqref{eq:3}.

Moreover,
\begin{equation}\label{eq:4}
  \mbox{$i\in L'(u)$ for every $u\in X_i$}.
\end{equation}
Otherwise, by \eqref{eq:3}, there are $u\in X_i$,  $w\in V(F)$ and  $v\in I[u,w]$ such that $\phi(w)=i$ and $i\notin L(v)$. Now, on the one hand, since $u,w\in X_i$ and $v\in I[u,w]$, the convexity of $X_i$ implies that $v\in X_i$. On the other hand, since $i\notin L(v)$, we know by \eqref{eq:3} that $v\notin X_i$. This contradiction proves \eqref{eq:4}.

Next, we now show that, for each $j\in[p]$,
\begin{equation}\label{eq:1}
\mbox{if $v'w'\in E$, with $w'\in X_j$ and $v'\notin X_j$, then
$j\notin L(v')$.}
\end{equation}
This is immediate if $v'w'\in F$, by the definition of $L(v')$.
Otherwise, there is $vw\in F$ such that $w\in X_j$, and $v, v'\in X_i$ for some $i\neq j$.
Lemma~\ref{l:frontierdistances} applied to the convex set $X_i$
and the edge $vw$ yields that  $d(v',v)<d(v',w)$; i.e., $v'\in X_{vw}$. Thus, the definition of $L(v')$ gives that $j\notin L(v')$, proving~\eqref{eq:1}.

We now prove \eqref{eq:2}. Consider $u\in X_i$ and $j\in[p]-\{i\}$. Let $w\in V(F)\cap X_j$ (as $G$ is connected, this set is non-empty) and let $P$ be a shortest path between $u$ and $w$. By construction, $P$ has some edge $vw'$ such that $v\notin X_j$ and $w'\in X_j$.  By~\eqref{eq:1}, we have that $j\notin L(v)$. As $v\in I[u,w]$, and as $\phi(w)=j$, the definition of $L'(u)$ implies that $j\notin L'(u)$. This completes the proof of~\eqref{eq:2}.

Therefore, a convex $p$-partition with skeleton $(F,\phi)$ exists if and only if the following conditions hold: (\emph{i}) $|L'(u)|=1$ for each vertex $u$ of $G$; (\emph{ii}) the parts of the corresponding partition are convex.
The time needed to find a convex $p$-partition with skeleton
$(F,\phi)$ is dominated by the time needed to compute the set of
distances $d_G(w,u)$ for each $w\in V(F)$ and each $u\in V$
which is linear for each set $F$.
Indeed, for each $w\in V(F)$ we can make a breadth first search starting at $w$
and delete $\phi(w)$ from $L(u)$, for  those $u$ in $X_{vw}$.
Similarly, we can construct lists $L'(u)$, by running for each $w\in V(F)$,
a breath first search starting at $w$ during which we delete $\phi(w)$ from $L(u)$
of all the descendent $u$ of a vertex $v$  for which $\phi(w)\notin L(v)$.
\qed\end{pf}

When given a connected bipartite graph $G$ and an integer $p$, we can decide
whether $G$ has a convex $p$-partition as follows. We first \emph{guess}
a \emph{candidate} skeleton $(F,\phi)$ with $\phi$  a $p$-coloring of $F$,
and then, by using Theorem~\ref{t:main}, we compute
in linear time the unique (if any) partition $\{X_1,\ldots,X_p\}$ associated to $(F,\phi)$. The choices for $(F,\phi)$ are bounded from above by a function
that depends only on $p$. In fact, if $(F,\phi)$ is a skeleton of some partition, then it must satisfy  the following properties.

\begin{itemize}
 \item The size of $F$ satisfies $|F|\in \{p-1,\ldots ,\binom{p}{2}\}$.
 \item The function $\phi$ is surjective.
 \item Identifying all vertices $v\in V(F)$ of the same colour under $\phi$ yields a connected simple graph.
\end{itemize}

From the first condition we know that there are at most $$\binom{\binom{n}{2}}{|F|}\leq \binom{n}{2}^{|F|}$$ choices    for the set $F$. Hence, there are
at most $\binom{p}{2}\binom{n}{2}^{\binom{p}{2}}$ choices in total.

From the second condition we know that there are roughly $\binom{p^2}{|F|} \leq p^{2p}$ functions $\phi$.
Since the problem of determining the convex $p$-partitions of a graph can be reduced in polynomial time to computing the convex $p'$-partitions of its components for $p'\in\{1,\ldots,p\}$ \cite{addsDM2011, addsENDM2007}, we conclude the following.

\begin{cor}\label{c:pconvexP}
For each fixed $p\geq 1$, all convex $p$-partitions of a bipartite graph
can be enumerated in polynomial time.
\end{cor}

\section*{Acknowledgements}

This work was conducted in the framework of the international cooperation project ``Operations Management and Operations Research: methodological problems and real-world applications'' funded by CONICET (Argentina) and CONICYT (Chile). L.N.~Grippo and M.D.~Safe were partially supported by UBACyT Grants 20020100100980 and 20020130100808BA, CO\-NICET PIP 112-200901-00178 and 112-201201-00450CO, and ANPCyT PICT-2012-1324 (Argentina). M.~Stein was supported by Fondecyt grant 1140766, and M.~Stein and
M.~Matamala  received support of  Fondo Basal PFB-03 and
N\'ucleo Milenio Information and Coordination in networks ICM/FIC P10-24F.






\end{document}